\documentstyle{amsppt}

\input amsppt1

\magnification=\magstep1

\NoRunningHeads

\pagewidth{15 truecm}

\pageheight{25 truecm}

\baselineskip=13pt

\loadbold

\nologo

\TagsOnRight

\parindent=20pt

\font\small=cmr9

\font\smallbf=cmbx9

\def\spb{\vskip 0.2truecm}

\def\mpb{\vskip 0.5truecm}

\def\bpb{\vskip 0.8truecm}

\hcorrection{0.5 truecm}

\NoBlackBoxes


\topspace{-1cm}

\mpb
\centerline{\bf Nonexistence of holomorphic submersions between complex unit balls}
\centerline{\bf equivariant with respect to a lattice} 
\centerline{\bf and their generalizations}

\mpb
\centerline{\bf Vincent Koziarz and Ngaiming Mok\footnote"*"{Research partially supported by the CERG grant HKU7034/04P of the the HKRGC, Hong Kong.}}

\bpb\bpb
{\narrower
\baselineskip=10pt
\parindent=15pt
\small\smallskip\noindent
{\smallbf Abstract.}
In this article we prove first of all the nonexistence of holomorphic 
submersions other than covering maps between compact quotients of complex unit 
balls, with a proof that works equally well in a more general equivariant 
setting.  For a non-equidimensional surjective holomorphic map between compact 
ball quotients, our method applies to show that the set of critical values must 
be nonempty and of codimension 1.  In the equivariant setting the line of 
arguments extend to holomorphic mappings of maximal rank into the complex
projective space or the complex Euclidean space, yielding in the
latter case a lower estimate on the dimension of the singular locus
of certain holomorphic maps defined by integrating holomorphic 1-forms. In 
another direction, we extend the nonexistence statement on holomorphic 
submersions to the case of ball quotients of finite volume, provided that the 
target complex unit ball is of dimension $m \ge$ 2, giving in particular a 
new proof that a local biholomorphism between
noncompact $m$-ball quotients of finite volume must be a covering map
whenever $m \ge$ 2. Finally, combining our results with Hermitian
metric rigidity, we show that any holomorphic submersion from a bounded
symmetric domain into a complex unit ball equivariant with respect to
a lattice must factor through a canonical projection to yield an
automorphism of the complex unit ball, provided that either the
lattice is cocompact or the ball is of dimension at least 2.
\smallskip}

\bpb
In the study of holomorphic mappings and rigidity problems on
compact quotients of bounded symmetric domains the case of 
$n$-ball quotients has always occupied a special place in 
terms of formulations of problems and methods of their solution.
(Here and in what follows by a quotient of a bounded symmetric 
domain we will always mean a quotient with respect to a
discrete torsion-free subgroup of biholomorphic automorphisms, and an 
$n$-ball quotient means a quotient of the $n$-dimensional complex 
unit ball $B^n$.) 
The method of harmonic maps of Siu ([Siu1], 1980) makes it possible to obtain 
holomorphic maps from harmonic maps into 
$m$-ball quotients under mild conditions, because the 
canonical K\"ahler-Einstein metric on the complex $n$-ball is of strictly
negative curvature in the dual sense of Nakano.  In the case 
where $m = 1$ any harmonic mapping $f: X \to C$ from a compact
K\"ahler manifold $X$ onto a compact Riemann surface $C$ of genus
$\geq 2$ of maximal rank at some point leads by the study of 
holomorphic foliations associated to $f$ to a holomorphic 
mapping $g: X \to C'$ onto some compact Riemann surface $C'$.
(Siu [Siu3]).  This means in general that representations of K\"ahler groups
into automorphism groups $\Bbb P\text{SU}(m,1)$ are associated to 
holomorphic objects.  When $X$ is an irreducible compact quotient
of a bounded symmetric domain, Margulis super-rigidity (Margulis
[Ma], 1977) or the method of Hermitian metric rigidity (Mok [Mok1], 1987)
implies that any holomorphic mapping $f: X \to Z$ into a compact
$m$-ball quotient $Z$ is necessarily constant, unless if $X$ is itself
an $n$-ball quotient for some $n$.  From this perspective 
holomorphic maps between compact quotients of 
complex unit balls are special and yet they are essential for 
completing our 
understanding of holomorphic maps between compact quotients 
of bounded symmetric domains, and more generally of linear 
representations of their fundamental groups.

\spb

For the case
of a holomorphic mapping $F: B^n \to B^m$, $n > m \geq 1$, between 
complex unit balls of maximal rank at some point and equivariant 
with respect to a representation of some cocompact lattice, it 
is generally believed that any such mapping must be singular
at certain points.  Especially there ought to be no holomorphic submersions
$f: X \to Z$ between compact quotients of complex unit balls 
$X := B^n/\Gamma$ and $Z := B^m/\Delta$, $n > m \geq 1$.  In this 
article for convenience we will refer to the latter as the 
Submersion Problem (for compact quotients of complex unit balls).  
The Submersion Problem for  $n = 2$ and  
$m = 1$ was settled by Liu ([Liu], 1996), in which more generally 
the nonexistence of regular holomorphic fibrations on compact 2-ball 
quotients was proven by means of Chern-number inequalities on 
surfaces arising from Teichm\"uller theory. In this article we resolve 
the Submersion Problem in all dimensions, proving more generally the 
nonexistence of holomorphic submersions from $B^n$ into $B^m$, 
$n > m \ge 1$, equivariant with respect to some representation
$\Phi: \Gamma \to \text{Aut}(B^m)$.

\spb

In Siu ([Siu2], 1984) there is a problem whether any holomorphic embedding 
between compact quotients
of complex unit balls must necessarily be totally geodesic provided 
that the domain manifold is of dimension $\geq 2$.  The problem can
be slightly generalized to allow for holomorphic immersions and 
the generalized problem will be referred to as the Immersion Problem  
(for compact quotients of the complex unit ball).  In Cao-Mok [CM]
the Immersion Problem was solved in the affirmative under 
the additional assumption that the complex dimensions $k$ 
resp. $\ell$ of the domain 
manifold resp. target manifold satisfy $\ell < 2k$.    
The starting point was an adaptation of an algebraic identity of 
Feder's ([Fe], 1965) on Chern classes (in the case of the projective space)
to the context of compact quotients of the complex unit ball.
Here Cao-Mok worked with the kernel of a closed nonnegative (1,1)-form 
arising from the second fundamental form of the holomorphic immersion. 
As it turns out, the same algebraic identity of Feder's can be adapted to 
the the study of holomorphic submersions from $B^n$ onto some open subset of $B^m$, 
$n > m \ge 1$.  Denoting by
$\omega_k$ the K\"ahler form of the complete K\"ahler manifolds
of constant holomorphic sectional curvature $-4\pi$ on the complex unit ball $B^k$, 
for any holomorphic mapping $F: B^n \to B^m$ by the Schwarz 
Lemma we have $F^*\omega_m \le \omega_n$, and the difference
$\omega_n - F^*\omega_m$ is a closed nonnegative (1,1) form which 
we will show to have a nontrivial kernel of dimension $m$ at each 
point.  We conclude that $F: B^n \to B^m$ is a Riemannian 
submersion and derive a contradiction from curvature properties
of the complex unit ball.

\spb
One of the new features that has come out from studying the 
Submersion Problem is a generalization of 
the solution of the Submersion Problem to subvarieties
of compact quotients of the complex unit ball.  As an illustration 
we prove that the set of regular values of a surjective holomorphic
map $f: X \to Z$ between compact quotients of complex unit balls 
does not contain any compact algebraic curve $C$ by studying  
the hypothetical holomorphic submersion on $f^{-1}(C) \subset X$. 
In relation to non-equidimensional surjective holomorphic maps 
$f: X \to Z$ between compact quotients of complex unit balls, it
should be noted that Siu raised in [Siu2] the problem whether 
such maps can exist when $\text{dim}(Z) \geq 2$.  There is at 
this point no convincing evidence one way or another.  Our results
are applicable to holomorphic maps from compact quotients of unit 
balls onto compact Riemann surfaces which were studied in Siu [Siu3]. 
At the same time, they provide some information on the critical 
values of such hypothetical maps when $\text{dim}(Z) \geq 2$ 
which could be useful for further investigation on the original 
question of Siu's hitherto unanswered.  Moreover, our methods are 
also applicable for surjective holomorphic maps from compact 
$n$-ball quotients into the complex projective space and compact complex tori (although 
the proofs are simpler), and in the latter case a formulation in 
the equivariant setting yields new information for critical values
of the integrals of linearly independent Abelian differentials on 
the covering complex unit ball.

\spb
In another direction our methods can be generalized to non-compact
complex hyperbolic space forms of finite volume.  (Here a complex
hyperbolic space form means a quotient of the complex unit ball by 
a torsion-free discrete group of automorphisms.) In relation to 
the Submersion Problem we show that any holomorphic submersion between 
non-compact complex hyperbolic space forms of finite volume must be 
a topological covering map (hence equidimensiona1) provided that the 
target manifold is of complex dimension $\ge 2$, noting that in 
the 1-dimensional case there are plenty of unramified holomorphic maps
between noncompact complex hyperbolic Riemann surfaces of finite volume 
which fail to be topological covering maps.  Our methods also lead
to a general structure theorem for holomorphic submersions from 
bounded symmetric domains to the complex unit ball $B^n$ equivariant 
with respect to a lattice $\Gamma$, showing that they must factor through 
a canonical projection to yield an automorphism of the complex 
unit ball itself, provided that either $\Gamma$ is uniform (i.e., cocompact) 
or $m \ge 2$ (and $\Gamma$ is a non-uniform lattice).

\spb\noindent
\smc{Acknowledgement.}
\rm
The current article is in part the fruit of research discussion 
between the authors while the first author was visiting the 
University of Hong Kong (HKU) in July 2007.  The latter would 
like to thank members and the secretarial staff of the Department 
of Mathematics of HKU for their hospitality.

\mpb

\noindent
{\bf \S1 Nonexistence of holomorphic submersions equivariant with 
respect to a cocompact lattice}

\noindent
We start by fixing some conventions and terminology.  A complex 
manifold is understood to be connected.  For complex manifolds $Y$ and $Z$, 
a mapping $h: Y \to Z$ is said to be a {\it holomorphic submersion\/} if and only 
if $h$ is holomorphic and $dh$ is of rank equal to $\text{dim}(Z)$ at 
every point.  A proper holomorphic submersion $h: Y \to Z$ (which is necessarily
surjective by the Proper Mapping Theorem) will be called a {\it regular holomorphic fibration\/} 
if and only if dim$(Z) \ge 1$ and fibers of $h$ are also
positive-dimensional.
When fibers of $h$ are connected, $h: Y \to Z$ realizes $Y$ 
as the total space of a regular family of compact complex manifolds 
$Y_t: = h^{-1}(t)$.  When one and hence any fiber of $h$ has $k > 1$
connected components, then there is a $k$-fold unramified cover $\rho: Z' \to Z$
and a regular holomorphic fibration $h': Y \to Z'$ with connected 
fibers such that $h = h' \circ \nu$. 

In the study of holomorphic mappings between compact quotients 
of bounded symmetric domains, in view of Hermitian metric 
rigidity [Mok1] what remains to be understood is essentially the case
where the domain manifold is a compact $n$-ball quotient.  In this 
direction there was the work of Cao-Mok [CM] on holomorphic immersions
between compact quotients of complex unit balls.
In our study of the Submersion Problem aiming at proving the 
nonexistence of regular holomorphic fibrations 
between compact quotients of complex unit balls it will be clear that 
the methods are equally applicable
to holomorphic submersions which are equivariant 
with respect to a 
cocompact lattice on the domain complex unit ball, and the discreteness
of the image of the underlying homomorphism does not play any essential 
role.  For this reason we will state and prove the first result in 
this broader context.  We have

\spb

\proclaim{Theorem 1}
For a positive integer $k$ denote by $B^k$ the $k$-dimensional complex unit
ball, and by {\rm Aut}$(B^k)$ its group of biholomorphic automorphisms.  
Let $\Gamma \subset$ {\rm Aut}$(B^n)$ be a cocompact lattice
of biholomorphic automorphisms.  Let $\Phi: \Gamma \to$ {\rm Aut}$(B^m)$
be a homomorphism and $F: B^n \to B^m$ be a holomorphic submersion
equivariant with respect to $\Phi$, i.e., $F(\gamma x) =  \Phi(\gamma)(F(x))$ for every $x \in B^n$ and for every $\gamma \in \Gamma$. 
Then,  $m=n$ and $F\in\hbox{\rm Aut}(B^n)$.   
\endproclaim

We observe also that torsion-freeness of
$\Gamma\subset\text{Aut}(B^n)$ is not needed in the hypothesis as one can always
pass to a torsion-free subgroup $\Gamma_0\subset\Gamma$ of finite index to
prove the theorem. As a corollary to Theorem 1 we have immediately

\proclaim{Corollary 1}
In the notations of Theorem 1 suppose furthermore that $\Gamma 
\subset$ {\rm Aut}$(B^n)$ is torsion-free, and write $X = B^n/\Gamma$.
Then, there does not exist any regular 
holomorphic fibration $f: X \to Z$
onto a compact $m$-ball quotient with $1\le m < n$.
\endproclaim

\spb

Corollary 1 for the case of $n = 2$ was established by Liu ([Liu], 1996) 
by means of geometric height inequalities in the case of curves over a 
complex function field of transcendence 1, i.e., for complex surfaces
fibered over a compact Riemann surface.  Liu's result was used 
by Kapovich [Ka]
to show that the compact 2-ball quotients constructed by Livn\'e [Liv]
have incoherent fundamental groups, i.e., they contain finitely
generated subgroups which are not finitely presentable.  The proof 
of Liu [Liu] makes use of inequalities arising from Teichm\"uller 
theory, and as such does not apply to the case of $\Gamma$-equivariant 
holomorphic maps into the unit disk.  It also does not generalize
to regular holomorphic fibrations of $n$-ball quotients over compact 
Riemann surfaces of genus $\geq$ 2 for $n \geq 3$, in which case
we are dealing with fibers of complex dimension $\geq 2$.

\spb

In 1965, Feder [Fe] proved that any holomorphic immersion 
$\tau: \Bbb P^k \to \Bbb P^{\ell}$ between complex projective 
spaces is necessarily a linear embedding whenever $\ell < 2k$.    
He did this by using Whitney's formula on Chern classes 
associated to the tangent sequence of the holomorphic sequence, 
thereby proving that the degree $\tau_*: H_2(\Bbb P^k,\Bbb Z) \to
H_2(\Bbb P^{\ell},\Bbb Z)$ must be 1 under the dimension restriction, a condition
which forces the vanishing of the $k$-th Chern class of the 
holomorphic normal bundle.  An adaptation of Feder's identity 
was used by Cao-Mok [CM] to study the Immersion Problem for the 
dual situation of holomorphic immersions between compact quotients 
of complex unit balls. For the Submersion Problem we have an 
associated short exact sequence, and the dual of that sequence 
is formally identical to the tangent sequence associated to 
holomorphic immersions, except that the role of the tangent
bundle is played by the cotangent bundle.  At the level of 
Chern classes Feder's identity remains applicable.
Representing first Chern classes in terms of the canonical 
K\"ahler-Einstein metrics and higher Chern classes by means 
of the Proportionality Principle of Hirzebruch, we are 
able to prove the following vanishing result on certain 
differential forms which serves as a starting point for the 
proof of Theorem 1.

\spb

Here and in what follows we will denote 
by $\omega_n$ the K\"ahler form of the complete 
K\"ahler-Einstein metric of constant holomorphic sectional 
curvature $-4\pi$ on the complex unit ball $B^n=\{z\in\Bbb C^n : \|z\|^2<1\}$. Writing the metric as $ds^2=2\hbox{\rm Re}\sum \, g_{i\bar j} dz^i\otimes d\bar z^j$ in local holomorphic coordinates, its curvature tensor is
$$R_{i\bar jk\bar \ell}=-2\pi(g_{i\bar j}\,g_{k\bar \ell}+g_{i\bar \ell}\,g_{k\bar j}).$$
The constant is chosen so that for the dual K\"ahler metric 
on the complex projective space $\Bbb P^n$, of constant holomorphic sectional 
curvature $4\pi$,  the K\"ahler form $\omega_{\Bbb P^n}$ 
represents the positive generator of $H^2(\Bbb P^n,\Bbb Z)$, as can be seen 
from the Gauss-Bonnet formula on $\Bbb P^1$.
For a quotient $X := B^n/\Gamma$ with respect to a torsion-free lattice 
$\Gamma$ we will write $\omega_X$ for the K\"ahler form induced
by $\omega_n$. When we have a holomorphic mapping $F: B^n \to B^m$
equivariant with respect to a representation of $\Gamma$, the 
(1,1)-form $F^*\omega_m$ is invariant under $\Gamma$, which 
induces a (1,1)-form on $X$ to be denoted by 
$\overline{\omega_m}$, bearing in mind that there is implicitly
the underlying holomorphic map $F$.  We have

\proclaim {Proposition 1}
Let $F: B^n \to B^m$, $X = B^n/\Gamma$ be as in the statement of Theorem 1.  
Then $\omega_n - F^*\omega_m$ is a nonnegative closed $(1,1)$-form 
on $X$, and $\big[\omega_X - \overline{\omega_m}\big]^{n-m+1} = 0$ 
as an $(n-m+1,n-m+1)$--cohomology class.  As a consequence 
$\big(\omega_X - \overline{\omega_m}\big)^{n-m+1} \equiv 0$ on $X$.
\endproclaim

\demo{Proof of Proposition 1}
By the choice of normalization of canonical metrics
it follows that the total Chern class of $\Bbb P^n$ is given by 
$\big(1+[\omega_{\Bbb P^n}]\big)^{n+1}$. 
By the Hirzebruch Proportionality Principle  
the total Chern class of the tangent bundle $T_X$ is given by 
$$
c(T_X) = \big(1-[\omega_X]\big)^{n+1}. \tag 1
$$
In particular $c_k(T_X)$ is a multiple of $[\omega_X]^k$ for 
$1 \leq k \leq n$.  From the $\Phi$-equivariant holomorphic submersion 
$F: B^n \to B^m$, the level sets of $F$ define a $\Gamma$-equivariant
holomorphic foliation which descends therefore to a holomorphic 
foliation $\Cal F$ on $X = B^n/\Gamma$.  We denote by $T_{\Cal F}$ the 
associated distribution on $X$.  Consider the short exact sequence
$0 \to T_{\Cal F} \to T_X \to N_{\Cal F} \to 0$ on $X$, which defines the 
holomorphic normal bundle $N_{\Cal F}$ of the foliation $\Cal F$.
Since $N_{\Cal F}$ is obtained by pulling back the tangent 
bundle of $B^m$ by $F: B^n \to B^m$ and descending to $X$, an analogue of (1)
applies to $N_{\Cal F}$, giving
$$
c(N_{\Cal F}) = \big(1 - \big[\overline{\omega_m}\big]\big)^{m+1}. \tag 2
$$
On the other hand, since the short exact sequence
is defined everywhere on $X$, by Whitney's formula we have  
$$
\big(1-[\omega_X]\big)^{n+1} = c(T_X) = c(T_{\Cal F})c(N_{\Cal F})
= c(T_{\Cal F})\big(1 - \big[\overline{\omega_m}\big]\big)^{m+1}. \tag 3
$$
Since $(B^n,\omega_n)$ and $(B^m,\omega_m)$ are both equipped with complete
K\"ahler-Einstein metrics of constant negative holomorphic sectional curvature
of the same negative constant, by the Schwarz Lemma we have 
$F^*\omega_m \leq \omega_n$. 
Proposition 1 is then a consequence of the following
elementary algebraic identity taken from Feder [Fe] for which we include a proof
for the sake of easy reference.
\enddemo

\proclaim{Lemma 1}
For the compact complex hyperbolic space form $X = B^n/\Gamma$ let 
$\alpha, \beta \in H^2(X,\Bbb R)$. Suppose for $1 \leq k \leq n-m$
there exists $\gamma_k \in H^{2k}(X,\Bbb R)$ such that
$(1+\alpha)^{n+1} = (1 + \gamma_1 + \cdots \gamma_{n-m})(1+\beta)^{m+1}$.
Then, $(\alpha-\beta)^{n-m+1} = 0$.
\endproclaim

\demo{Proof}
Let $\gamma \in \bigoplus_{k=0}^{n} H^{2k}(X,\Bbb R)$ be the formal 
quotient $(1+\alpha)^{n+1}(1+\beta)^{-(m+1)}$.  Writing $\gamma_k$ for 
the component of degree $2k$ in $\gamma$, the notation is consistent 
with those in the statement of Lemma 1, and we have $\gamma_{n-m+1} = 0$.
We compute $\gamma_{n-m+1}$ formally.  For an element $\delta$ of the 
graded cohomology groups of even degrees we write $\delta_k$ for its
element of degree $2k$. We have
$$
\gather
0 = \gamma_{n-m+1} = \big((1+\alpha)^{n+1}(1+\beta)^{-(m+1)}\big)_{n-m+1} \\
= \sum_{k+\ell =n-m+1} (-1)^{\ell}\frac{(n+1)!}{k!\ (n-k+1)!}\frac{(m+\ell)!}{m! \ \ell !}
\alpha^k\beta^{\ell}\\   
= \sum_{k+\ell =n-m+1} (-1)^{\ell} \frac{(n+1)!}{k!\ (n-k+1)!}\frac{(n-k+1)!}{m! \ \ell !}
\alpha^k\beta^{\ell}\\
= \sum_{k+\ell =n-m+1} (-1)^{\ell}\frac{(n+1)!}{(n-m+1)!\ m!}
\frac{(n-m+1)!}{k! \ \ell !} \alpha^k\beta^{\ell} \\
= \frac{(n+1)!}{(n-m+1)!\ m!}(\alpha-\beta)^{n-m+1} \\
\endgather
$$
as desired. \quad \quad $\square$
\enddemo

We are now ready to give a proof of Theorem 1.

\demo{Proof of Theorem 1}
Consider the closed (1,1)-form $\rho := \omega_X -
\overline{\omega_m}$ on $X$.  By Proposition 1, $\rho \geq 0$
and $\rho^{n-m+1} = 0$.  Since by definition $F^*\omega_m$ vanishes
on the level set $F^{-1}(w)$ for any $w \in B^m$, on $B^n$ the 
(1,1)-form $\omega_n - F^*\omega_m$ agrees with $\omega_n$ on each 
level set of $F$, of dimension $n-m$, so that $\rho$ must have 
at least $n-m$ positive eigenvalues everywhere on $X$.  The 
identity $\rho^{n-m+1} \equiv 0$ implies that at every point
$x \in X$, all other eigenvalues of $\rho(x)$ are zero.  In 
other words, there exists an $m$-dimensional complex vector 
subspace $H_x \subset T_{X,x}$ transversal to $T_{\Cal F, x}$ such 
that $\rho(x)|_{H_x} \equiv 0$.

\spb
From the short exact sequence $0 \to T_{\Cal F} \to T_X \to
N_{\Cal F}\to 0$ there are two different ways to endow the holomorphic vector 
bundle $N_{\Cal F}$ with a Hermitian metric.  First, endowing the holomorphic
tangent bundle $T_X$ with the Hermitian metric $g_X$ defined by the K\"ahler form 
$\omega_X$,  $N_{\Cal F} = T_X/T_{\Cal F}$ inherits a Hermitian metric $h$ as a Hermitian 
holomorphic quotient vector bundle.  On the other hand, $\Gamma$ acts on the Hermitian
holomorphic vector bundle $F^*T_{B^m}$, and the latter descends to a Hermitian 
holomorphic vector bundle on $X$ which is isomorphic to $N_{\Cal F}$ as a 
holomorphic vector bundle.  From this we obtain another Hermitian metric $h'$ on $X$.
We argue that the two Hermitian metrics $h$ and $h'$ agree with each other.
To see this at $x \in X$ write $T_{X,x} = T_{\Cal F, x} \oplus H_x$ as a complex
vector space.  When measured against the Hermitian inner product $g_X(x)$ on $T_{X,x}$, the
direct summands $T_{\Cal F, x}$ resp. $H_x$ are eigenspaces
of the Hermitian inner product defined by $\rho(x)$ corresponding to the 
eigenvalues 1 resp. 0.  As a consequence they must be orthogonal to each other.
For an element $\eta \in N_{\Cal F,x}$, the norm $\|\eta\|_h$ of $\eta$ with respect to $h$
is given by the minimum norm of $\|\widetilde{\eta}\|_{g_X}$ with respect to $g_X$, as 
$\widetilde{\eta}$ ranges over all (1,0)-vectors at $x$ which projects to $\eta$
modulo $T_{\Cal F,x}$.  Now that $T_{X,x} = T_{\Cal F, x} \oplus H_x$ is an orthogonal 
decomposition, $\eta$ lifts to $\eta_0 \in H_x$, and $\|\eta\|_h$ is nothing other
than $\|\eta_0\|_{g_X}$.  However, since $\rho|_{H_x} \equiv 0$, we have $\overline{\omega_m}|_{H_x}
\equiv \omega_X|_{H_x}$, so that $\|\eta_0\|_g = \|\eta_0\|_{h'}$, proving that
$h \equiv h'$. In other words, $F:B^n\to B^m$ is an isometric submersion in the sense of
Riemannian geometry.

\spb

Therefore, if $m=n$, $F$ induces a (local) isometry of $B^n$ onto itsfelf, so that $F$ sends local, hence global, geodesics of $B^n$ to geodesics of $B^n$. It must be injective and proper, hence a biholomorphism.

\spb

We suppose now that $n>m$ and we wish to get a contradiction. We can compute the curvature tensor of the Chern connection on $N_{\Cal F}$ associated to the metric $h=h'$ in two different ways.  Denote by $\sigma\in C^\infty_{1,0} (X,\hbox{\rm Hom}(T_{\Cal F},N_{\Cal F}))$ the second fundamental form of the Hermitian holomorphic vector subbundle $T_{\Cal F} \subset T_X$ with respect to the K\"ahler-Einstein metric $g_X$ and by $\sigma^*\in C^\infty_{0,1} (X,\hbox{\rm Hom}(N_{\Cal F},T_{\Cal F}))$ its adjoint. Let $x\in X$ and let $(\xi_1,\dots,\xi_n)$ be an orthonormal basis of $T_{X,x}$ with respect to $g_X(x)$ such that $\xi_{m+1},\dots,\xi_n$ belong to $T_{{\Cal F},x}$. Writing $\Theta$ resp. $\Theta'$ for the curvature tensor of $(T_X,g_X)$ resp. $(N_{\Cal F},h)$ at $x$ we have by a classical computation of Griffiths
$$\Theta'_{j\bar k\lambda\bar\mu}=\Theta_{j\bar k\lambda\bar\mu}+\bigl\langle\sigma^*_{\bar\xi_k}(\xi_\lambda),\overline{\sigma^*_{\bar\xi_j}(\xi_\mu)}\bigr\rangle_{g_X}$$
for any $j,k\in\{1,\dots,n\}$ and any $\lambda,\mu\in\{1,\dots,m\}$, where we have identified $N_{{\Cal F},x}$ with $H_x=T_{{\Cal F}_x}^\perp$ (see for example [De]). But recall that
$$\Theta_{j\bar k\lambda\bar\mu}=-2\pi(\delta_{j\bar k}\delta_{\lambda\bar\mu}+\delta_{j\bar\mu}\delta_{\lambda\bar k})$$
and since $h=h'$, $\Theta'$ is the pullback by $dF$ of the curvature tensor of $T_{B^m}$:
$$\eqalign{\Theta'_{j\bar k\lambda\bar\mu}& = -2\pi(\delta_{j\bar k}\delta_{\lambda\bar\mu}+\delta_{j\bar\mu}\delta_{\lambda\bar k}),\ \ \forall\, j,k,\lambda,\mu\in\{1,\dots,m\}\cr
&=  0 \ \ \hbox{\rm in other cases.}\cr}$$
Comparing the different relations, we easily deduce that the family 
$$\bigg \{\frac{\sigma^\star_{\bar\xi_j}(\xi_\lambda)}{\sqrt{2\pi}}\,:\,(j,\lambda)
\in \{m+1,\dots,n\}\times \{1,\dots,m\}\bigg \}$$
of vectors of $T_{{\Cal F},x}$ is orthonormal.
If $m\geq 2$, it is clear (because of the dimensions) that the latter property can never be verified and we get the desired contradiction.
If $m=1$ then $N_{\Cal F}$ is a line bundle. For any $x\in X$, if $\xi,\xi'\in T_{{\Cal F},x}$ and $\eta\in N_{{\Cal F},x}$
$$\bigl\langle\sigma_{\xi}(\xi'),\overline\eta\bigr\rangle_{h}=\bigl\langle\xi',\overline{\sigma^*_{\bar\xi}(\eta)}\bigr\rangle_{g_X}$$
and the above property implies that the restriction of $\sigma: T_{\Cal F}\otimes T_{\Cal F}\longrightarrow N_{\Cal F}$ is everywhere non-degenerate. In other words, if $\Omega_{\Cal F}$ denotes the dual bundle of $T_{\Cal F}$, then $T_{\Cal F}$ and $\Omega_{\Cal F}\otimes N_{\Cal F}$ are isomorphic as smooth complex vector bundles (actually, they are isomorphic as holomorphic vector bundles since $\sigma$ is in fact holomorphic, but we will not need that in what follows). In particular, their first Chern classes coincide, hence $c_1(T_{\Cal F})=-c_1(T_{\Cal F})+(n-1)(-2\,\overline{\omega_1})$, i.e., $c_1(T_{\Cal F})=-(n-1)\,\overline{\omega_1}$. Now, still as cohomology classes,
$$-(n+1)\,\omega_X=c_1(T_X)=c_1(T_{\Cal F})+c_1(N_{\Cal F})=-(n+1)\,\overline{\omega_1}
$$
and this is impossible because $F$ is a holomorphic Riemannian submersion with $n>m$.
The proof of Theorem 1 is complete.\quad \quad $\square$
\enddemo

\mpb\noindent
{\bf \S2 On the singular loci of holomorphic submersions}

\noindent
From the proof of Theorem 1 we obtain also the following result regarding regular
holomorphic fibrations
which also applies to the case where the target manifold is a 
compact complex torus or the complex projective space.  The 
complex unit ball equipped with a complete K\"ahler-Einstein metric 
is sometimes referred to as the complex hyperbolic space.  In this 
context a quotient of the complex unit ball by a torsion-free 
discrete group of biholomorphic automorphisms is sometimes called
a complex hyperbolic space form.

\proclaim{Theorem 2}
Let $n > m \geq 1$. Let $\Gamma \subset$ {\rm Aut}$(B^n)$ be a torsion-free cocompact lattice
of biholomorphic automorphisms, $X: = B^n/\Gamma$.  Let $Z$ be an $m$-dimensional 
compact complex hyperbolic space form, compact complex torus or complex projective space.
Let $f: X \to Z$ be a surjective holomorphic map and denote by $E \subset Z$ the 
smallest subvariety such that $f$ is a regular holomorphic fibration over $Z-E$.  Then,
there is no compact analytic subvariety of positive dimension in $Z-E$.  In particular,
$E \subset Z$ is of complex codimension 1.
\endproclaim

\demo{Proof}
On $f^{-1}(Z-E)$ we have the short exact sequence of holomorphic vector bundles 
$0 \to T_f \to T_X \to N \to 0$, where $N$ denotes the holomorphic normal bundle 
of the holomorphic foliation $\Cal F$ defined by the relative tangent bundle 
$T_f$, so that $c(T_X) = c(T_f)c(N)$ holds on it.
Let $Q \subset Z-E$ be an irreducible complex-analytic curve.  Restrict the short exact sequence 
to the compact complex-analytic subvariety
$f^{-1}(Q)$, even if $Q$ may have singularities. Denote by $\omega_Z$   
the closed (1,1)-form on $Z$ such that $-(m+1)\omega_Z$ is the first Chern form of a canonical 
K\"ahler-Einstein metric $g_Z$ on $Z$.  This is consistent with notations in the proof of Proposition 1
if $Z$ is a compact quotient of the $m$-ball, and is defined for the case of compact complex tori
and the complex projective space in such a way that Proposition 1 remains applicable.
We conclude from $c_{n-m+1}(N) = 0$ and Proposition 1
that $[\omega_X - f^* \omega_Z]^{n-m+1} = 0$ as a cohomology class on $f^{-1}(Q)$.

\spb
In case $Z$ is a compact complex torus or the complex projective space, this is already a 
contradiction since $\omega_Z$ is nonpositive and hence $\omega_X - f^* \omega_Z$ 
is strictly positive, and the subtle case is the complex hyperbolic case, where $\omega_X - f^*\omega_Z$
is only known to be nonnegative. In this case, arguing as in the proof of Theorem~1, we find that $f_{|f^{-1}(Q)}:f^{-1}(Q)\rightarrow Q$  is a Riemannian submersion above each nonsingular point of $Q$, if $f^{-1}(Q)$ (resp. $Q$) is endowed with the metric induced by $g_X$ (resp. $g_Z$). Let $z\in Q$ be a regular point of $Q$ and consider the exact sequence of bundles on the fiber $P=f^{-1}(z)$
$$0\longrightarrow T_P\longrightarrow T_{f^{-1}(Q)_{|P}}\longrightarrow N_P \longrightarrow 0$$
where $T_P$ is the tangent bundle to the manifold $P$ and $T_{f^{-1}(Q)}$ is the tangent bundle to (the regular part of) $f^{-1}(Q)$, both equipped with the metrics induced by $g_X$. We also endow the normal bundle $N_P= T_{f^{-1}(Q)_{|P}}/T_P$ with the quotient metric denoted by $h$. Let $x\in P$ and let $(\xi_1,\dots,\xi_{n-m+1})$ be an orthonormal basis of $T_{f^{-1}(Q),x}$
such that $\xi_2,\dots,\xi_{n-m+1}$ are tangent to the fiber $P$. Writing $\Phi$ resp. $\Phi'$ for the curvature tensor of $T_{f^{-1}(Q)_{|P}}$ resp. $N_P$, we have for any $2\leq j,k\leq n-m+1$,
$$0=\Phi'_{j\bar k1\bar 1}=\Phi_{j\bar k1\bar 1} +\bigl\langle{\tau}^*_{\bar\xi_k}(\xi_1),\overline{{\tau}^*_{\bar\xi_j}(\xi_1)}\bigr\rangle_{g_{X}},
$$
where the second fundamental form $\tau\in C^\infty(P,\hbox{\rm Hom}(S^2T_P,N_P))$ of the submanifold $P\subset X$ is considered as an element of $C^\infty_{1,0}(P,\hbox{\rm Hom}(T_{P},N_P))$, and $\tau^*\in C^\infty_{0,1} (P,\hbox{\rm Hom}(N_{P},T_P))$ is its adjoint (again, we identify $N_P$ with $T_P^\perp$).
In particular, because of the Chern-Weil formula ($T_{f^{-1}(Q)_{|P}}$ being a holomorphic subbundle of $T_{X_{|P}}$), for any $\eta\in N_{P,x}$ and any $\xi\in T_{P,x}$,
$$2\pi\|\xi\|^2\|\eta\|^2\leq-\Phi_{\xi\bar\xi \eta\bar\eta}= \bigl\langle{\tau}^*_{\bar\xi}(\eta),\overline{{\tau}^*_{\bar\xi}(\eta)}\bigr\rangle_{g_X}=\bigl\langle\tau_{\xi}({\tau}^*_{\bar\xi}(\eta)),\overline\eta\bigr\rangle_h.
$$
Since the line bundle $N_P$ is trivial along $P$, $\tau$ can be seen as a symmetric bilinear form on $T_P$. 
By the previous inequality, $\tau$ is non-degenerate on $T_{P,x}$ and this is true for any $x\in P$. As a consequence, $T_{P}$ and $T_P^*$ are isomorphic as smooth bundles and therefore $c_1(T_P)=0$.
But $c_1(T_P)<0$ because $X$ is K\"ahler-Einstein with negative Einstein constant and $T_P$ is a holomorphic subbundle of ${T_X}_{|P}$, so we obtain a contradiction, proving Theorem~2.\quad \quad $\square$
\enddemo

\spb

From the statement of Theorem 2 we deduce readily

\proclaim{Corollary 2}
A compact complex hyperbolic space form does not admit any regular holomorphic 
fibration over a compact K\"ahler  
manifold of constant holomorphic sectional curvature \text{\rm (}i.e., a compact hyperbolic 
space form, a compact complex torus, or a complex projective space\text{\rm )}.  As a consequence,
a compact complex hyperbolic space form does not admit any regular holomorphic 
fibration over a compact Riemann surface.
\endproclaim

\noindent
{\smc Remarks.}
Corollary 2 follows already from an easy extension of the 
proof of Theorem 1 
to cover the case where the target manifold is the complex 
Euclidean space or the complex projective space (as included
in the proof of Theorem 2). 
The last statement of Corollary 2 was covered 
by Liu [Liu] in the special case when the domain manifold is
of dimension 2.  Corollary 2 leaves 
open the question whether nontrivial regular holomorphic fibrations over higher 
dimensional base manifolds can exist on     
compact complex hyperbolic space forms of dimension $\ge 3$.

\spb
We include some results which follow readily from modifications of the proof of Theorem 1. In the proof of Theorem 1, where we derive a 
contradiction by assuming that the $\Phi$-equivariant holomorphic mapping
$F: B^n \to B^m$ is a holomorphic submersion, the argument actually 
works to arrive at a contradiction provided that the set of singularities 
of $F$, i.e., the subset $\text{Sing}(F) \subset B^n$ over which $dF$ 
is not of maximal rank, is sufficiently small.  More precisely, if 
$\text{Sing}(F)$ is of dimension $< m-1$, then removing $\text{Sing}(F)$ has 
no effect on $c_{n-m+1}(N_{\Cal F})$.  This is so because a generic hyperplane 
section of $X$ of dimension $n-m+1$ does not intersect the locus $S \subset X$,
where $S := \text{Sing}(F)/\Gamma$.  From this we deduce

\proclaim{Theorem 3}
Let $n > m \geq 1$ and let $\Gamma \subset$ {\rm Aut}$(B^n)$ be a cocompact lattice
of biholomorphic automorphisms.  Let $\Phi: \Gamma \to$ {\rm Aut}$(B^m)$
be a homomorphism into the automorphism group {\rm Aut}$(B^m)$ of $B^m$, 
and $F: B^n \to B^m$ be a nonconstant holomorphic
map equivariant with respect to $\Phi$ which is a holomorphic submersion
at some point.  Let $\text{\rm Sing}(F) \subset B^n$ 
be the $\Gamma$-invariant subset consisting of points where $F$ fails to
be a submersion, {\rm i.e.}, where $dF$ is of rank $< m$, which descends
to a complex-analytic subvariety $S \subset X$.  Then, $S$ is nonempty
and $\text{\rm dim}(S) \geq m-1$.   
\endproclaim

\noindent
{\smc Remarks.} 
Note that Theorem 3, applied to the special case of holomorphic 
submersions $f: X \to Z$ between compact complex hyperbolic space forms, does not 
imply Theorem 2.  In fact, from the statement of Theorem 3 it does 
not even follow that the image of $f(S) = E \subset Z$ is of dimension $m-1$, i.e., 
that $E$ is of codimension 1 on $Z$, since we do not know that the fiber of 
$f|_S: S \to E$ over a general point of $E$ is isolated.  Beyond saying that 
$E \subset Z$ is of codimension 1, Theorem 2 actually suggests that the 
codimension-1 components of $E$ resemble an ample divisor.

\spb
For the study of $\Gamma$-equivariant holomorphic submersions $F:B^n\to B^m$, the arguments remain valid 
(with a simpler proof) when the target manifold is replaced by the Euclidean space or the complex 
projective space (cf. the proof of Theorem~2). In particular, 
suppose $\nu_1, \cdots \nu_m$ are $m$ linearly independent holomorphic 1-forms
on a compact quotient $X = B^n/\Gamma$ of the complex unit ball $B^n$, an analogue 
of Theorem 3 remains valid for the holomorphic mapping $F: B^n \to \Bbb C^m$ obtained by 
integrating pull-backs of the holomorphic 1-forms $\nu_1. \cdots \nu_m$ by the 
universal covering map $\pi: B^n \to X$.  In other words, we have for the mapping 
$F$ obtained as integrals of Abelian differentials the following statement on 
singularities of the meromorphic foliation defined by level sets of $F$.

\proclaim{Theorem 4}
Let $n > m \geq 1$. Let $\Gamma \subset$ {\rm Aut}$(B^n)$ be a cocompact lattice
of biholomorphic automorphisms, $X: = B^n/\Gamma$. Let $\nu_1,\cdots,\nu_m$ be 
$m$ holomorphic 1-forms on $X$ which are linearly independent at a general point 
of $X$.  Let $S \subset X$ be the subvariety at which $\nu_1,\cdots,\nu_m$ fail 
to be linearly independent.  Then, $\text{\rm dim}(S) \geq m-1$. 
\endproclaim

\mpb\noindent
{\bf \S3 Generalization to complex hyperbolic space forms of finite volume}

\noindent
In this section, we prove a version of Theorem~1 in the case where $\Gamma\subset$ {\rm Aut}$(B^n)$ is a non-uniform lattice. This means that the quotient $B^n/\Gamma$ is non-compact, but the volume of $X$ 
with respect to $\omega_X$ is finite. Our arguments are quite elementary, in the sense that they do not make use of any compactification of $X$.

\proclaim{Theorem 1'}  
Let $\Gamma \subset$ {\rm Aut}$(B^n)$ be a lattice
of biholomorphic automorphisms.  Let $\Phi: \Gamma \to$ {\rm Aut}$(B^m)$
be a homomorphism and $F: B^n \to B^m$ be a holomorphic submersion
equivariant with respect to $\Phi$. Suppose either $m\geq 2$ or $\Gamma
\subset${\rm Aut}$(B^n)$ is cocompact.  
Then,  $m=n$ and $F\in\hbox{\rm Aut}(B^n)$.   
\endproclaim

\spb
When $\Gamma \subset$ {\rm Aut}$(B^n)$ is a non-uniform lattice, 
it is necessary to impose the condition $m \ge 2$.
In fact, for $m = 1$ there are plenty of unramified holomorphic
maps between Kobayashi-hyperbolic punctured Riemann surfaces 
which are not topological coverings (cf. Remarks after the proof).

\demo{Proof of Theorem 1'} We only have to prove the theorem when $\Gamma \subset$ {\rm Aut}$(B^n)$ 
is torsion-free and $X := B^n/\Gamma$ is non-compact. In this situation, we can argue exactly as in the proof of Theorem~1 if we show that
Proposition~1 is still valid. In fact, the proof of the vanishing of $\big[\omega_X - \overline{\omega_m}\big]^{n-m+1}$ 
as an $(n-m+1,n-m+1)$--cohomology class goes along the same line because in particular the Hirzebruch Proportionality Principle remains valid, but we need to explain why $\big(\omega_X - \overline{\omega_m}\big)^{n-m+1} \equiv 0$ on $X$.
Let us remark that when $m=1$, $n-m+1=n$, and in this case the vanishing of the class above does not give any information since $H^{2n}(X,{\Bbb R})=0$ if $X$ is non-compact.

In order to continue the proof, we need to recall some facts about the geometry of the manifold $X$ (see W.~M.~Goldman's book~[Go] or [KM] for more details). It is the union of a compact part and of a finite number of disjoint cusps. Each cusp $C$ is diffeomorphic to a product $N\times[0,+\infty)$ where $N$ is a compact quotient of a horosphere $HS$ of $B^n$ centered at a point $\infty\in\partial B^n$. The fundamental group $\Gamma_C$ of $C$ may be identified with the stabilizer in $\Gamma$ of the horosphere $HS$.

Let $(z,v,t)\in{\Bbb C}^{n-1}\times{\Bbb R}\times{\Bbb R}$ be horospherical coordinates associated to $HS$. The 1-form $\varsigma=\frac{1}{2\pi}\,e^{-t}(2\hbox{\rm Im}\langle\langle z,dz\rangle\rangle-dv)$, as well as $dz$, $t$ and $dt^2$ are invariant by $\Gamma_C$ and the metric $g_X$ takes the form
$$g_X=\frac{1}{2\pi}\,(dt^2+\varsigma^2+4e^{-t}\langle\langle dz,dz\rangle\rangle)$$
in the cusp $C$. The fundamental remark is that on the cusp $C$, $d\varsigma=\omega_X$ (the invariance of $\varsigma$ by $\Gamma_C$ allows to go down on $C$) and, because of the form of the metric, $|\varsigma|_{g_X}$ is constant.

Let us go back to the proof of the proposition. Let $K$ be a compact subset of $X$ which contains in its interior the compact part of $X$. Let $\chi\in C^\infty(X,{\Bbb R})$ be a smooth function vanishing on the compact part of $X$, and equal to 1 on $X\backslash K$. The 1-form $\alpha=\chi\varsigma$ is well-defined on $X$ (the definition of $\varsigma$ of course depends on the cusp). Moreover, the 2-form $(\omega_X-d\alpha)$ has compact support in $X$. Therefore,
$$\int_X\bigl(\omega_X - \overline{\omega_m}\bigr)^{n-m+1}\wedge\omega_X^{m-2}\wedge (\omega_X-d\alpha)=0$$
since $\bigl(\omega_X - \overline{\omega_m}\bigr)^{n-m+1}$ vanishes in $H^{n-m+1,n-m+1}(M,{\Bbb R})$ and since it is integrated against a $d$-closed form with compact support.

As $(X,g_X)$ is complete, there exists an exhaustive sequence $(K_\nu)_{\nu\in{\Bbb N}}$ of compact subsets of $X$ and 
smooth cut-off functions $\psi_\nu :X\to [0,1]$ which are identically equal to one on $K_\nu$, vanish on $X\backslash K_{\nu+1}$, and verify $|d\psi_{\nu}|_{g_X}\leq 2^{-\nu}$.
Now, by the Schwarz Lemma, $\bigl|\bigl(\omega_X - \overline{\omega_m}\bigr)^{n-m+1}\bigr|_{g_X}$ is uniformly bounded by some constant. Noting that $(X,\omega_X)$ is of finite volume,
$$\lim_{\nu\rightarrow+\infty}\int_X \bigl(\omega_X -\overline{\omega_m} \bigr)^{n-m+1}\wedge\omega_X^{m-2}\wedge d \psi_\nu\wedge\alpha=0$$
and then
$$\int_X \bigl(\omega_X - \overline{\omega_m}\bigr)^{n-m+1}\wedge\omega_X^{m-2}\wedge d\alpha=\lim_{\nu\rightarrow+\infty}\int_X \bigl(\omega_X - \overline{\omega_m}\bigr)^{n-m+1}\wedge\omega_X^{m-2}\wedge d( \psi_\nu\,\alpha)=0.$$
We immediately deduce that
$$
\int_X\bigl(\omega_X - \overline{\omega_m}\bigr)^{n-m+1}\wedge\omega_X^{m-1}=0
$$
and then, that the conclusion of Proposition~1 is true. The proof of Theorem 1' is complete.\quad \quad $\square$
\enddemo

\spb\noindent
{\smc Remarks.}

\noindent 1) Theorem 1', as it is stated, is trivially false when $X$ is non-compact and $n=1$. For example, let $\Upsilon\subset\hbox{\rm Aut}(B^1)$ be any cocompact lattice, $Y= B^1/\Upsilon$. Let $X$ be the same Riemann surface as $Y$ with $p>0$ punctures. Then, there exists a lattice $\Gamma\subset\hbox{\rm Aut}(B^1)$ such that $X$ is biholomorphic to $B^1/\Gamma$ but the embedding $X\hookrightarrow Y$, which is also a submersion, is not isometric.

\spb
\noindent 2) Other than standard examples, very few examples of representations of lattices of 
$\hbox{\rm Aut}(B^n)$ into $\hbox{\rm Aut}(B^m)$ are known.
Nevertheless, one can find in Deligne-Mostow [DM] some examples --- based on a construction of R.~Livn\'e [Liv] --- of non-trivial holomorphic maps $f:X\longrightarrow Y$ between compact complex hyperbolic manifolds, with $n=2$ and $m=1$ (cf. also Kapovich [Ka]). 
It is conceivable that the method of construction as expounded in [DM] 
can also lead to examples with $n=3$ and $m=1$, although no such examples are available in the literature.
We also refer to the examples of Mostow in the case $n=m=2$ detailed by Toledo~[To].

\spb
\noindent 3) When $\Gamma$ is arithmetic, the Satake-Borel-Baily compactification $\overline{X}$, which is
projective-algebraic, is obtained by adding a finite
number of cusps, which are isolated singularities
of $\overline{X}$.  The proof of Theorem 1' can
then be obtained by slicing $\overline{X}$ to
obtain hyperplane sections which avoid the cusps.
When $m \ge 2$ then $n-m+1 \le n-1$.  The vanishing
of $(\omega_X - \overline{\omega_m})^{n-m+1}$ on
all such hyperplane sections implies its identical
vanishing on $X$, which gives Theorem 1'.
The same argument applies in the case of non-arithmetic quotients.  Here it is known that $X$ can be compactified by adding a finite number
of points (Siu-Yau [SY]), but the  proof of its
projective-algebraicity does not seem to be
available in the literature.  In place of overloading
the article with writing down a self-contained proof
of the latter, we have chosen to present the more
elementary argument here.

\spb
\noindent 4)
For $\Gamma \subset \text{Aut}(B^n)$ a non-uniform lattice, the case of Theorem 1'
for $n = m \ge 2$ is already non-trivial.  It implies in particular that a local 
biholomorphism from a non-compact complex hyperbolic space form of finite volume 
into a complex hyperbolic space form is necessarily a covering map, which 
was established in Mok [Mok2, p.174ff.] in a much more elaborate way 
by the method of Hermitian metric rigidity applied to certain homogeneous 
holomorphic vector bundles.

\mpb\noindent
{\bf \S4 Structure of holomorphic submersions from finite volume
quotients of bounded symmetric domains into the complex
unit ball}

\noindent In this section we consider the general structure of
holomorphic submersions of quotients of bounded
symmetric domains $\Omega$ into the complex unit ball.
Let $\Omega=\Omega_1\times\dots\times\Omega_q$ be the decomposition
of $\Omega$ into a product of irreducible bounded
symmetric domains $\Omega_i, 1 \le i \le q$. We assume that each $\Omega_i$ is noncompact and denote by $\hbox{\rm Aut}(\Omega_i)$ the group of biholomorphic automorphisms of $\Omega_i$.
Let $\Gamma\subset\hbox{\rm Aut}(\Omega)$ be a lattice. After passing to a subgroup of $\Gamma$ of finite index, one can always assume that $\Gamma$ is torsion-free and belongs to $\hbox{\rm Aut}_0(\Omega)$, the identity component of $\hbox{\rm Aut}(\Omega)$. Then, there exists a partition $I_1,\dots,I_p$ of $\{1,\dots,q\}$ and irreducible lattices $\Gamma_{I_k}\subset \Pi_{i\in I_k}\hbox{\rm Aut}_0(\Omega_{i})$ such that $\Gamma=\Pi_{i=1}^p \Gamma_{I_k}$.
The fact that the $\Gamma_{I_k}$ are irreducible means that for any proper subset $J$ of $I_k$ the projection of $\Gamma_{I_k}$ into $\Pi_{j\in J}\hbox{\rm Aut}_0(\Omega_j)$ is dense (see [Ra, Cor. 5.21]).
We shall denote by $X$ the quotient manifold $\Omega/\Gamma$. Note that the tangent bundle of $X$ has a natural decomposition $T_X=T_{X,1}\oplus\dots\oplus T_{X,q}$ coming from the decomposition of $\Omega$.

The following result is a consequence of Theorems~1 and 1', and of Hermitian metric rigidity (see [Mok1] to which we will frequently refer below).

\proclaim{Theorem 5}  
Let $\Omega$, $\Gamma$ and $X$ be as above. Let $\Phi:\Gamma\to\hbox{\rm Aut}(B^m)$
be a homomorphism and $F: \Omega \to B^m$ be a holomorphic submersion equivariant with respect to $\Phi$. 
Suppose that $m\geq 2$ or that $X$ is compact. Then, there exists a Hermitian symmetric space $\Omega'$ such that $\Omega=B^m\times\Omega'$ and $F$ is the natural projection onto $B^m$ composed with an element of $\hbox{\rm Aut}(B^m)$.  
\endproclaim

\demo{Proof} We show first that there exists $\ell\in\{1,\dots,q\}$ such that the restriction of $dF$ to $\oplus_{j\neq \ell}T_{X,j}$ vanishes identically. Let $z^{(i)}=(z^{(i)}_1,\dots,z^{(i)}_{n_i})$ be Euclidean coordinates on $\Omega_i$ in terms of its Harish-Chandra embedding and $z=(z^{(1)},\dots,z^{(q)})$ be the corresponding coordinates on $\Omega$. We endow $\Omega_i$ with the (unique up to some constant) Bergman metric $h_i$ and denote by $\pi_i:\Omega\to\Omega_i$ the natural projection. Consider now the K\"ahler metric $h=\sum\pi_i^*h_i+F^*g_m$ on $\Omega$. We note once for all that, because of the Schwarz Lemma, $h$ is dominated by a constant multiple of $\sum\pi_i^*h_i$. Indeed, the holomorphic sectional curvature of $h_i$ is negative and bounded from below, and the holomorphic sectional curvature on $B^m$ is negative and constant.

The metric $h$ goes down on $X$. From the proof of Theorem~4 in [Mok1] it follows easily that the 2-tensor $F^*g_m$ is given by
$$
F^*g_m=2\,\hbox{\rm Re}\,\Bigl(\sum_{{1\leq i\leq q}\atop{1\leq j,k\leq n_i}} a^{(i)}_{j\bar k}(z^{(i)}) \,dz^{(i)}_j\otimes d\bar z^{(i)}_k
+2\sum_{{{1\leq i<i'\leq q}\atop{1\leq j\leq n_i}}\atop{1\leq k\leq n_{i'}}}b^{(i,i')}_{j\bar k}(z)\,dz^{(i)}_j\otimes d\bar z^{(i')}_k
\Bigr)
$$
for some functions $a^{(i)}_{j\bar k}:\Omega_i\to{\Bbb C}$,  $b^{(i,i')}_{j\bar k}:\Omega\to{\Bbb C}$. In fact, if $\Omega_i$ has at least rank two then the matrix of functions $(a^{(i)}_{j\bar k})$ defines the canonical K\"ahler-Einstein metric on $\Omega_i$ (the (ir)reducibility of $\Gamma$ does not play any role in that case). If $\Omega_i$ is of rank one, the functions $a^{(i)}_{j\bar k}$ only depend on the coordinates of $\Omega_i$, and they define the canonical K\"ahler-Einstein metric on $\Omega_i$ whenever the cardinality of the multi-index $I_k$ containing $i$ is at least 2. Since the $(1,1)$-form associated to $F^*g_m$ is $d$-closed, $b^{(i,i')}_{j\bar k}$ must be holomorphic in the $z^{(i)}$-coordinates and antiholomorphic in the $z^{(i')}$-coordinates.

Suppose now that there exists $x\in\Omega$, two integers $i\neq i'$ and tangent vectors  $\xi=\partial/\partial z^{(i)}_j,\eta=\partial/\partial z^{(i')}_{j'}\in T_{\Omega,x}$ whose images by $d_xF$ do not vanish. We endow the vector bundle $T_{X,i}\oplus T_{X,i'}$ with the metric $\pi^* h_i+\pi^* h_{i'}+(F^*{g_m})_{|T_{X,i}\oplus T_{X,i'}}$ and denote by $\Theta$ its curvature tensor. In other words, $T_{X,i}\oplus T_{X,i'}$ is identified with a holomorphic vector subbundle of $T_{X,i}\oplus T_{X,i'}\oplus F^*T_{B^m}$ by the embedding
$i(\eta_i,\eta_j) = (\eta_i,\eta_j,df(\eta_i+\eta_j))$,
from which it inherits a Hermitian metric. Because of the curvature decreasing property of the curvature for holomorphic subbundles, and since bisectional curvatures of $B^m$ are strictly negative, we have $\Theta_{\xi\bar\xi\eta\bar\eta}\neq0$ (see the proof of Lemma 2 in [Mok1]). But in the previous notations,
$$\gather
(F^*{g_m})_{|T_{X,i}\oplus T_{X,i'}}=2\,\hbox{\rm Re}\,\Bigl(\sum_{1\leq j,k\leq n_i} a^{(i)}_{j\bar k}(z^{(i)}) \,dz^{(i)}_j\otimes d\bar z^{(i)}_k+\\
\sum_{1\leq j,k\leq n_{i'}} a^{(i')}_{j\bar k}(z^{(i')}) \,dz^{(i')}_j\otimes d\bar z^{(i')}_k
+2\sum_{{1\leq j\leq n_i}\atop{1\leq k\leq n_{i'}}}b^{(i,i')}_{j\bar k}(z)\,dz^{(i)}_j\otimes d\bar z^{(i')}_k
\Bigr)
\endgather
$$
and from the properties of the functions $a^{(i)}_{j\bar k}$, $a^{(i')}_{j\bar k}$ and $b^{(i,i')}_{j\bar k}$, we deduce that the value of $\Theta_{\xi\bar\xi\eta\bar\eta}$ is not affected by the presence of the term $F^*g_m$ and hence that $\Theta_{\xi\bar\xi\eta\bar\eta}=0$. This is a contradiction,
thus $dF$ must vanish in the direction of all the $\Omega_i$ but one at each point of $\Omega$. By holomorphicity of $F$, $dF$ must vanish identically in the direction of each factor except one, say $\Omega_1$. One can therefore regard $F$ as a map from $\Omega_1$ to $B^m$.

\spb

As a consequence, we can assume that $\Gamma=\Gamma_{I_1}$ with $1\in I_1$. If $I_1=\{1\}$ then Theorem~1 of [Mok1] applied to $(\Omega_1/\Gamma_1,h_1)$ and the metric $h_1+F^*g_m$ implies that $\Omega_1$ must be a complex unit ball, otherwise $X_1 = \Omega_1/\Gamma_1$ is of rank $\ge 2$, and $F^*g_m = c\,h_1$ for some constant $c$
by Hermitian metric rigidity, which is impossible as $g_m$ is of strictly negative bisectional curvature.
So we can apply Theorem~1 or 1', and Theorem 5 is proved in this case. In other cases, it follows directly from Theorem~4 of [Mok1] that $F^*g_m=c\,h_1$ for some global constant $c$. Then, necessarily, $\Omega_1=B^m$ and $F:B^m\to B^m$ must be injective and proper, so $F\in\hbox{\rm Aut}(B^m)$.\quad \quad $\square$
\enddemo

\spb\noindent
{\smc Remarks.}
\noindent If, in the previous theorem, the image of $\Phi$ is supposed to be discrete (which is the case whenever $F$ is induced by a holomorphic map from $X$ to a complex hyperbolic space form) then $\Gamma$ must be reducible. More precisely, $\Gamma=\Gamma_0\times\Gamma'$ where $\Gamma_0\subset\hbox{\rm Aut}(B^m)$ and $\Gamma'\subset\hbox{\rm Aut}_0(\Omega')$ are lattices. Indeed, if this were not the case, the projection $\hbox{\rm pr}_1(\Gamma)$ into $\hbox{\rm Aut}(B^m)$ would be dense. This is impossible since $\hbox{\rm pr}_1(\Gamma)$ is conjugate to $\Phi(\Gamma)$ in $\hbox{\rm Aut}(B^m)$, by Theorem~5.

\bpb
\noindent
{\bf References}

\parindent=30pt

\item{[CM]}
Cao, H.-D.; Mok, N. \quad
Holomorphic immersions between compact
hyperbolic space forms,
{\it Invent. Math.}
{\bf 100} (1990), 49-61.

\item{[De]}
Demailly, J.-P. \quad
Complex analytic and algebraic geometry,
``Open content book'' available on the webpage of J.-P. Demailly.

\item{[DM]}
Deligne, P.; Mostow, G.D. \quad
{\it Commensurabilities among Lattices in $PU(1,n)$\/},
Princeton University Press, Princeton, New Jersey, 1993.

\item{[Fe]}
Feder, S. \quad
Immersions and embeddings in complex projective spaces,
{\it To\-pology}
{\bf 4} (1965), 143-158.

\item{[Go]}
Goldman, W. M. \quad
Complex hyperbolic geometry,
The Clarendon Press, Oxford University Press, 1999

\item{[Ka]}
Kapovich, M. \quad
On normal subgroups in the fundamental groups 
of complex surfaces, at 
{\it http://www.math.utah.edu/~kapovich/EPR/koda.pdf}.
Preprint 1998.

\item{[KM]}
Koziarz V.; Maubon J. \quad
Harmonic maps and representations of non-uni\-form lattices of $\hbox{\rm PU}(m,1)$,
{\it Ann. Inst. Fourier (Grenoble)}
{\bf 58} (2008), 507-558.

\item{[Liu]}
Liu, K. \quad
Geometric height inequalities,
{\it Math. Res. Lett.}
{\bf 3} (1996), 693-702.

\item{[Liv]}
Livn\'e, R. \quad
On certain covers of the universal elliptic curve,
Ph. D. Thesis, Harvard University, 1981

\item{[Ma]} 
Margulis, G. A. \quad
Discrete groups of motion of manifolds of nonpositive curvature,
{\it A.M.S. Transl.}(2)
{\bf 109} (1977), 33-45.

\item{[Mok1]}
Mok, N. \quad
Uniqueness theorem of Hermitian metrics of seminegative curvature
on quotients of bounded symmetric domains,
{\it Ann. Math.}
{\bf 125} (1987), 105-152.

\item{[Mok2]}
Mok, N. \quad
{\it Metric Rigidity Theorems on Hermitian Locally Symmetric Manifolds\/},
Series in Pure Mathematics Vol.{\bf 6}, 
World Scientific, Singapore-New Jersey-London-Hong Kong, 1989.

\item{[Siu1]}
Siu, Y.-T. \quad
The complex analyticity of harmonic maps and the strong rigidity 
of compact K\" ahler manifolds,
{\it Ann. Math.}
{\bf 112} (1980), 73-111.

\item{[Siu2]}
Siu, Y.-T. \quad
Some recent results in complex manifold theory related to vanishing 
theorems for the semipositive case, {\it Arbeitstagung Bonn} 1984,
Lect. Notes Math., Vol. 1111, Springer-Verlag 1985, 
Berlin-Heidelberg-New York, pp.169-192.

\item{[Siu3]}
Siu, Y.-T. \quad
Strong rigidity for K\"ahler manifolds and the construction of 
bounded holomorphic functions, in 
{\it Discrete Groups in Geometry and Analysis}, papers in honor of 
G. D. Mostow, ed. R. Howe, Birkh\"auser, Boston-Basel-Stuttgart,
1987, pp.124-151.

\item{[SY]}
Siu, Y.-T. and Yau, S.-T. \quad
Compactification of negatively curved complete K\"ahler
manifolds of finite volume, in {\it Seminar on
Differential  Geometry},  pp. 363--380, Ann. of Math.
Stud., 102, Princeton Univ. Press, Princeton, N.J., 1982.

\item{[To]}
Toledo, D. \quad
Maps between complex hyperbolic surfaces,
{\it Geom. Dedicata}
{\bf 97} (2003), 115-128.

\mpb\noindent
Vincent Koziarz,
Institut Elie Cartan, Universit\'e Henri Poincar\'e,
B. P. 239, F-54506 Vand\oe uvre-l\`es-Nancy Cedex, France
\newline
(E-mail : koziarz\@iecn.u-nancy.fr)

\spb\noindent
Ngaiming Mok, The University of Hong Kong, Pokfulam Road, Hong Kong
\newline
(E-mail: nmok\@hkucc.hku.hk)

\bye